\documentclass[12pt,a4paper]{article}
\usepackage{hyperref}
\usepackage[german,english]{babel} 
\usepackage{amsthm}
\usepackage{amsmath}
\usepackage{amsfonts}
\usepackage{amssymb}
\usepackage{enumerate}
\usepackage{amscd}
\usepackage{graphicx}
\usepackage{makeidx}
\usepackage[all]{xy}

\makeindex

\newcounter{Def}[section]

\theoremstyle{break}
\newtheorem{Proposition}[Def]{Proposition}

\newtheorem{Remark}[Def]{Remark}
\newtheorem{Lemma}[Def]{Lemma}
\newtheorem{Theorem}[Def] {Theorem}
\newtheorem{Kor}[Def] {Corollary}

\newenvironment{Proof}[1][\quad]
{\mbox{}\\\noindent \textbf{Proof:{\hspace{0.7cm}
#1}\\}}{\hfill$\Box$\\}

\newenvironment{Definition}{ \refstepcounter{Def}\mbox{}\\ 
\noindent\sl\textbf{Definition
\arabic{section}.\arabic{Def}} }
{\vspace{0.3cm}}

\def\C          {\mathbb C}

\def\be                 {\begin{equation}}
\def\ee                 {\end{equation}}

\def\id               {{\rm id}}

\def\tn                 {\tilde{n}}

\def\tg                 {\tilde{g}}

\def\tx                 {\tilde{x}}
\def\ty                 {\tilde{y}}

\def\tk                 {\tilde{k}}
\def\id                 {\mathrm{id}}
\def\coker                {\mathrm{coker}}
\def\im                 {\mathrm{Im}}
\def\inx                {\in \mathcal{X}_2}

\def\gm                 {\mathcal{X} = (\mathcal{X}_1 , \mathcal{X}_2, \mu,         \partial)}
\def\x                  {\mathcal{X}}
\def\bx                 {\overline{\mathcal{X}}}

\def\irr                {\Lambda_{\mathcal{M}(\mathcal{X})}}
\def\vk                 {\underline{\mathbf{0}}}
\def\teins              {\underline{\mathbf{1}}}
\def\c                  {\mathcal{C}}   
\def\hom                {\mathrm{Hom}}
\def\mod                {\mathrm{Mod}}
\def\m                  {m}

\def\ind                {\mathrm{Ind}_{\vk}}
\def\mx                 {\mathcal{M}(\mathcal{X})}
\def\bmx                {\overline{\mx}}
\def\mbx                {\mathcal{M}(\overline{\mathcal{X}})}
\def\msx                {\mathcal{M}(\mathcal{X}')}

\def\kca                {K^* \! \rtimes_{\hat \mu} \! C}
\def\K                                  {\frac{1}{|K|}}

\begin{document}


\thispagestyle{empty}
\begin{flushright}
   {\sf ZMP-HH/10-4}\\
   {\sf Hamburger$\;$Beitr\"age$\;$zur$\;$Mathematik$\;$Nr.$\;$363}\\[2mm]
   March 2010
\end{flushright}
\vskip 2.0em
\begin{center}\Large 
MODULAR CATEGORIES FROM FINITE CROSSED MODULES
\end{center}\vskip 1.4em
\begin{center}
  Jennifer Maier and 
  ~Christoph Schweigert\,\footnote{\scriptsize 
  ~Email addresses: \\
  $~$\hspace*{2.4em}Jennifer.Maier@uni-hamburg.de, 
  Christoph.Schweigert@uni-hamburg.de}
\end{center}

\vskip 3mm

\begin{center}\it
  Organisationseinheit Mathematik, \ Universit\"at Hamburg\\
  Bereich Algebra und Zahlentheorie\\
  Bundesstra\ss e 55, \ D\,--\,20\,146\, Hamburg
\end{center}
\vskip 2.5em
\begin{abstract} \noindent
It is known that finite crossed modules provide premodular tensor
categories. These categories are in fact modularizable. We construct the
modularization and show that it is equivalent to the module category of a 
finite Drinfeld double.
\end{abstract}

\setcounter{footnote}{0} \def\thefootnote{\arabic{footnote}} 


\renewcommand{\labelenumi}{(\roman{enumi})}

\section{Introduction}

Modular tensor categories and, more generally, premodular
tensor categories arise as representation categories
of certain (weak) Hopf algebras, certain nets of von Neumann
algebras and suitable classes of vertex algebras.
They have found numerous applications, including
the construction of invariants of three-manifolds and links,
the construction of low-dimensional quantum field theories
and the construction of gates in topological quantum computing. The simplest algebraic
object whose representation category is a premodular tensor
category is a finite crossed module.

\begin{Definition}\mbox{} \\
A finite \emph{crossed module} consists of two finite
groups $\x_1$ and $\x_2$, together with a (right)action
$\mu$ of $\x_1$ on $\x_2$ by group automorphisms,
written as $\mu(m,g)=m^g$, and
a group homomorphism, called the {\em boundary map},
$\partial:\ \x_2\to\x_1$ that satisfies
$$\ \partial(m^g)= g^{-1}(\partial m) g \,\, 
\text{ and }\,\, m^{\partial n}= n^{-1}mn 
\qquad \text{ for all }\,\, m,n\in\x_2 \,\text{ and } \,g\in\x_1.$$
\end{Definition}

It follows immediately from the definition that the kernel of $\partial$ is 
a central subgroup of $\x_2$\ and that the image of
$\partial$ is a normal subgroup of $\x_1$. This definition reduces to
the usual definition of the Drinfeld double $\mathcal{D}(G)$
of a finite group $G$ if the boundary map $\partial$ is the 
identity and the action $\mu$ is given by conjugation.

Our results hold over any algebraically closed field $k$ of characteristic zero.
Any finite crossed module $\x$ gives rise to a category
of representations over  $k$ which we denote 
by $\mx$.
The objects $(V,P,Q)$ of the category $\mx$ are finite-dimensional $\x_2$-graded $k$-vector spaces\\ $V=\oplus_{m\in\x_2} V_m$
with an action $Q:\x_1\to \mathrm{Aut}(V)$ of $\x_1$ such that
$$ P(m) Q(g)= Q(g) P(m^g) \text{ for all } m\in\x_2 \text{ and } g\in\x_1 \,\, . $$
Here $P(m)$ is the projection to the $m$ graded component.
Morphisms are required to preserve the $\x_2$-grading and the 
$\x_1$-action. In other words, we consider the category of
$\x_1$-equivariant vector bundles on $\x_2$ of finite rank.

We endow the category $\mx$ with the structure of a tensor
category.
The tensor product  is the usual tensor product of vector spaces, where the grading on $V \otimes W$ is given by $(V \otimes W)_m = \bigoplus_{nl= m} V_n \otimes W_l$ and the action of $\x_1$ on $V \otimes W$ is the diagonal action.  
The boundary map $\partial$
gives the additional structure of a braided tensor category
on $\mx$. Braiding isomorphisms are given by 
\begin{align}
R_{VW}:V \otimes W & \rightarrow W \otimes V
\label{Verzopfung}\\
v \otimes w & \mapsto \sum_{m \in \mathcal{X}_2} Q_W(\partial m)w \otimes P_V(m)v\nonumber \,\, .
\end{align}
Bantay has shown \cite{bantay}, that
together with the dualities inherited from the category 
of finite-dimensional $k$-vector spaces, where on the dual space the grading is defined by $(V^*)_m = (V_{m^{-1}})^*$ and the action is given by $Q^*(g) = Q(g^{-1})^*$,
the category $\mx$ has the structure of a premodular tensor category.

\begin{Definition} \mbox{} \\[-1.8em]
\begin{enumerate}
\item
Let $k$ be an algebraically closed field of characteristic zero.
A \emph{premodular tensor category} over $k$ is an abelian, $k$-linear, semi-simple ribbon category $\mathcal C$
such that
\renewcommand{\labelenumi}{\arabic{enumi})}
\begin{enumerate}
\item The tensor product is linear in each variable and the tensor unit is absolutely
simple, $\mathrm{End}(\mathbf{1}) = k$.
\item There are only finitely many isomorphism classes of simple objects, indexed by a set $\Lambda_{\c}$.
\end{enumerate} 
\renewcommand{\labelenumi}{(\roman{enumi})}
\item
The braiding on $\mathcal C$ allows to define the  
\emph{S-matrix} with entries in the field $k$
\begin{align}
s_{XY} := tr(R_{YX} \circ R_{XY}) \,\, ,
\label{smat}
\end{align} 
where $X,Y \in \Lambda_{\c}$.
A premodular category is called \emph{modular}, if the S-matrix is invertible.
\end{enumerate}
\end{Definition} 

We refer to \cite{BK,Ka} for the notion of a ribbon category.
For a detailed discussion of the premodular tensor category
$\mx$, including its character theory, we refer to \cite{bantay}.

 Modular tensor categories
are of particular interest, since they allow the construction of a topological field theory
\cite{RT,Tu} and thus of invariants of three-manifolds and of knots and links. The category $\mx$ associated to a crossed module is known to be modular, iff the boundary
map $\partial$ is an isomorphism \cite[Proposition 5.6]{N}.
In this case, it is equivalent to the representation category
of the Drinfeld double of a finite group. Equivalent categories
appear as representation categories of holomorphic
orbifold theories, see e.g.\ \cite{DVVV}.

Brugui\`eres \cite{br} (see also \cite{M1}) has introduced the notion of
modularization that associates to any premodular tensor category 
(obeying certain conditions) a modular tensor category. This
tensor category is unique up to equivalence of braided
tensor categories. The categories associated to crossed modules obey these conditions \cite{bantay}; hence
the question arises whether crossed modules provide a source of new modular tensor categories.
A first main result of this note is a negative answer to this question in Theorem \ref{Haupt}: the modularization
yields a modular tensor category equivalent to the category for the Drinfeld double
of $\x_2/\mathrm{ker}\,\partial \cong \mathrm{Im}\,\partial$.

Brugui\`eres has also given an explicit modularization procedure which is  based on a Tannakian subcategory
of the premodular tensor category. As a second main result of this note, we determine in Proposition \ref{Gruppendarst}
the group corresponding to the Tannakian subcategory to be a semi-direct product
$$ G(\x):=(\mathrm{ker}\,\partial)^{*} \! \rtimes_{\hat \mu} \! (\mathrm{coker}\, \partial) \,\, .$$
Here $(\mathrm{ker}\,\partial)^{*}$ is the group of characters
of the finite abelian group
$\mathrm{ker}\,\partial$; the semidirect product is 
explained in equation (\ref{G(X)}).
 The regular representation
of $G$ then provides a commutative special symmetric Frobenius algebra  
in the premodular
tensor category $\mx$. By general arguments, the category of left modules over this algebra
is a modular tensor category, see Proposition \ref{Induction}.

This note is organized as follows. In section 2 we recollect a few more aspects of crossed modules
and their representation category from \cite{bantay} and describe explicitly the full Tannakian subcategory
of transparent objects.
The transparent object corresponding to the regular representation of $G(\x)$ is shown in
Section 3 to be a commutative special symmetric Frobenius algebra $\vk$.
We describe the modularization  functor as the induction functor with respect to the algebra $\vk$.
In section 4, this description is used to construct an explicit equivalence of categories from the
modularization to the representation category of a Drinfeld double.

\section{Premodular categories from finite crossed modules}

We start by summarizing some more aspects of the
premodular category $\mx$ defined in the previous section. For any
object $V\in\mx$, the character is defined as the function
$$\begin{array}{rll}
\psi:\ \x_2\times \x_1 &\to& k \\
(m,g)&\mapsto &\ \mathrm{tr}_V (P(m)\ Q(g))\ .
\end{array} $$
The character theory for finite crossed modules largely parallels
(and in fact generalizes) the character theory of finite
groups \cite{bantay}. In particular, Maschke's theorem and
orthogonality relations hold:
for a general field $k$, the irreducible characters are orthogonal for
the non-degenerate symmetric bilinear form
\[ <\psi_1, \psi_2> := \frac{1}{|\x_1|} \sum_{g \in \x_1, m \inx} {\psi_1(m,g^{-1})} \psi_2(m,g) \,\, . \] 
For $k=\C$, the irreducible characters are orthonormal with respect
to the 
 hermitian scalar product
\[ (\psi_1, \psi_2) := \frac{1}{|\x_1|} \sum_{g \in \x_1, m \inx} \overline{\psi_1(m,g)} \psi_2(m,g) \,\, . \] 

We introduce two particularly important objects in $\mx$. To this end, we introduce the notation 
$K := \ker \,\partial, C := \coker \,\partial$ and
$I:=\im\,\partial$.
\begin{enumerate}
\item The tensor unit $\teins$ is defined on a one-dimensional
$k$-vector space in the graded component for $1\in\x_2$ and 
with trivial action of $\x_1$.

\item The vacuum object is the triple
$\vk = (V_{\vk}, P_{\vk},       Q_{\vk})$ 
with the $k$-vector space 
\[V_{\vk} = k[\mathrm{ker}\,\partial] 
\otimes k(\mathrm{coker}\,\partial) \equiv  k[K] 
\otimes k(C)\,\, .\]  

On the distinguished
basis $(x \otimes \delta_{Iy})_{x \in K, Iy \in C}$ we set
for $m \inx, \; g \in \x_1$:
\begin{eqnarray*}
P_{\vk}(m) (x \otimes \delta_{Iy}) = \delta(m^y,x)(x \otimes \delta_{Iy}) & \quad Q_{\vk}(g) (x \otimes \delta_{Iy}) = (x \otimes \delta_{Igy})
\end{eqnarray*} 

\end{enumerate}

A direct calculation using the explicit form (\ref{Verzopfung})
of the braiding gives the S-matrix defined in (\ref{smat}) in terms of the characters:
the entry corresponding to the irreducible representations $p,q \in \irr$ is
\[s_{pq} = \sum_{m,n \inx} \psi_p(m,\partial n) \psi_q(n, \partial m)
\,\, . \]
It is convenient to introduce a normalization factor to obtain the matrix:
\[S_{pq} := \frac{s_{pq}}{|\x|} = \frac{1}{|\x|} \sum_{m,n \inx} \psi_p(m,\partial n) \psi_q(n, \partial m)\]
with $|\x| := |\x_1|\cdot|\mathrm{ker}\,\partial|=|\x_2|\cdot |\mathrm{coker}\,\partial|$.
For later reference, we associate to each $p\in\Lambda_{\mx}$
the number
$$\ \omega_p:=\frac1{d_p}\sum_{m\in\x_2} \psi_p(m,\partial m)
\,\, ,\ $$
where $d_p$ is the categorical dimension of 
$p$\ (which coincides with the dimension of the underlying
vector space). It gives the eigenvalue of the twist $\theta_p$ on the simple object $p$, and it can be shown to satisfy the equality
\begin{align}
\psi_p(m,g\partial m) = \omega_p \cdot \psi_p(m,g)
\,\,.\label{15}
\end{align}

\begin{Remark}\mbox{}\\[-1.8em]
\begin{enumerate}
\item
Given any simple object $p \in \irr$, we have 
$S_{\teins p} = \frac{d_p}{|\x|}$, where $\teins$ is the
tensor unit.
\item
The multiplicity $\mu_p = \dim_{k} \hom(p,\vk)$ of the irreducible representation
$p$ in  $\vk$ equals
\begin{align}
\mu_p = D[S^2]_{\teins p}\,\, ,\label{25}
\end{align}
where $D := |\mathrm{coker}\,\partial|\cdot|\mathrm{ker}\,\partial|$. 
This follows in a straightforward calculation by expressing
the multiplicity in terms of characters as 
$\mu_p = (\psi_p, \psi_{\vk})$ and then using orthogonality
relations to compare with the matrix element of $S^2$.
\end{enumerate}
\end{Remark}

We recall the following definition from \cite{br}:

\begin{Definition} \\
An object $X$\ of a braided tensor category $\c$ is called
\emph{transparent}, if the equation $R_{Y,X} = R^{-1}_{X,Y}$
holds for all $Y\in\c$. We denote the set of
isomorphism classes of simple transparent objects by
$T_{\c}$.
\end{Definition}

The following observations are straightforward:

\begin{Remark}\label{2.3}\mbox{}\\[-1.8em]
\begin{enumerate}
\item
Direct summands of transparent objects and direct sums
of transparent objects are transparent.
\item
The vacuum object $\vk$ of $\mx$ is transparent.
This follows from a straightforward calculation using the
explicit form (\ref{Verzopfung}) of the braiding.
\end{enumerate}
\end{Remark}

\begin{Lemma}\label{vielfinvk}\mbox{} \\
An irreducible representation  $p \in \irr$ is 
a direct summand of the vacuum object $\vk$, if and only if the $p$-th row of the
S-matrix is collinear to the row  $(S_{\teins q})_{q \in \irr}$.

In this case, the multiplicity $\mu_p$ 
equals 
$\alpha = \mu_p = d_p$, where
$S_{pq} = \alpha S_{\teins q}$. Moreover, the twist
on any such  simple object with $\mu_p > 0$ is the identity.
\end{Lemma}

\begin{Proof}
Suppose, there is an $\alpha \in k$ such that 
$S_{pq} = \alpha S_{\teins q}$ for all $q \in \irr$. 
Specializing to $q=\teins$ yields
$\frac{d_p}{|\x|} =  \frac{\alpha}{|\x|}$ and hence the first
identity $\alpha = d_p > 0$. We compute the
multiplicity $\mu_p$:
\begin{align*}
\mu_p & \stackrel{(\scriptsize\ref{25})}{=} D[S^2]_{\teins p}
= D \sum_{q \in \irr} S_{\teins q} S_{p q}
= \alpha D \sum_{q \in \irr} S_{\teins q}^2
= \alpha \frac{D}{|\x|^2} \sum_{q \in \irr} 
 d_p^2=\ \alpha\,\, \,,
\end{align*}
where in the last step we have used a generalization
of Burnside's Theorem \cite{bantay} for the characters.

Conversely, suppose $\mu_p > 0$ so that $p$ is
a direct summand of the transparent object $\vk$ and hence,
by Remark \ref{2.3},
transparent itself. For any $q\in \irr$ we conclude
$R_{qp} \circ R_{pq} = \id_{p \otimes q}$ and thus
\begin{align*}
S_{pq} = \frac{1}{|\x|} \mathrm{tr}(\id_{p \otimes q}) = \frac{d_p d_q}{|\x|} = d_p S_{\teins q}  \,\, .
\end{align*} 
The equalities 
\begin{align*}
\mu_p & = \frac{1}{|\x_2|} \sum_{n \inx, m \in K} \psi_p(m,\partial n)
 \stackrel{(\scriptsize\ref{15})}{=} \frac{1}{|\x_2|} \sum_{n \inx, m \in K} \psi_p(m,\partial nm^{-1}) \, \omega_p\\[.2em]
& = \frac{1}{|\x_2|} \sum_{\tn \inx, m \in K} \psi_p(m,\partial \tn) \, \omega_p = \mu_p \omega_p
\end{align*}
show that $\mu_p > 0$ implies  $\omega_p = 1$.                
\end{Proof}

For a premodular category $\c$, consider the set of
isomorphism classes of those simple objects $X\in\c$
for which the row $(s_{XY})_{Y \in \Lambda_{\c}}$
is collinear with the row $(s_{\teins Y})_{Y \in \Lambda_{\c}}$ of the tensor unit:
\begin{align*}
M_{\c} = \{X \in \Lambda_{\c}\, \, |\,\, \forall \, Y \in \Lambda_{\c}:  s_{XY} = \dim X \dim Y\}  \,\, .
\end{align*}

\begin{Kor}\label{mt}\mbox{} \\
We have the following identities for the category $\mx$:
\begin{enumerate}
\item $M_{\mx} = T_{\mx}$.
\item $\theta_X = \id_X$ for all transparent objects $X$.
\item $\sum_{p \in T_{\mx}} (d_p)^2 = |\mathrm{ker}\,\partial||\mathrm{coker}\,\partial|$ with $d_p = \dim p$.
\end{enumerate}
\end{Kor}

\begin{Proof}
\vspace{-0.8cm}
\begin{enumerate}
\item According to Lemma \ref{vielfinvk}, any simple object
$p \in M_{\mx}$ is contained in the transparent object
$\vk$ and thus transparent itself. The other inclusion is
obvious. 
\item Lemma \ref{vielfinvk} and the first assertion of this
corollary imply $\theta_p = \id_p$ 
for all $p \in T_{\mx}$. The assertion follows since a
transparent object is a direct sum of simple transparent 
objects. 
\item The definition of $\mu_p$ and Lemma \ref{vielfinvk} 
imply 
\begin{align*}
|\ker\,\partial|\cdot|\coker\,\partial| & = \dim\vk
 = \sum_{p \in \irr} \mu_p d_p
 = \sum_{p \in T_{\mx}} (d_p)^2 \,\, .
\end{align*}
\end{enumerate}
\end{Proof}

Brugui\`eres' modularity criterion \cite[Proposition 1.1]{br}
asserts that a premodular category $\c$ is modular if and only if 
$M_{\c} = \{\mathbf{1}\}$. As an application we obtain:

\begin{Proposition}\label{modular}\mbox{} \\
The category $\mx$ is modular, if and only if the boundary map
$\partial$ is a bijection. In this case, $\mx$ is equivalent to the 
representation category of a Drinfeld double.
\end{Proposition}

\begin{Proof}
Lemma \ref{vielfinvk} implies that the row 
$(S_{pq})_{q \in \irr}$ is collinear with 
$(S_{\teins q})_{q \in \irr}$, if and only if $p$ has non-vanishing
multiplicity in $\vk$. For the tensor unit, we have
multiplicity $\mu_{\teins} = d_{\teins} = 1$.

If the boundary map $\partial$ is a bijection, we have
$\vk\cong\teins$ and, according to Brugui\`eres' criterion, 
the category $\mx$ is modular. If $\partial$ is not
bijective, we have $\dim V_{\vk} > 1$ and $\vk$ contains at 
least one simple object that is not isomorphic to the tensor 
unit $\teins$. Brugui\`eres' criterion now implies 
that the category is not modular. 
\end{Proof}

For a proof of this assertion that does not directly use
Bruigi\`eres' criterion, we refer to \cite[Proposition 5.6]{N}.

We will now explain why the premodular category $\mx$\ is 
modularizable \cite{bantay}.
To this end, we repeat some definitions of \cite{br}:

\begin{Definition} \mbox{} \\[-1.8em]
\begin{enumerate}
\item
An object $X$ of a category $\c$ is called a \emph{retract} of an
object $Y\in\c$, if there are morphisms $\iota:X \rightarrow Y$ 
and $\pi: Y \rightarrow X$ such that $\pi \circ\iota = \id_X$. 

\item
A functor $F:\c \rightarrow \c'$ is called \emph{dominant}, if for
every object $X\in\c'$ there exists an object $Y\in\c$ such that
$X$ is a retract of $F(Y)$.

\item
A \emph{modularization} of a premodular category $\c$ is a dominant
ribbon functor $F:\mathcal{C} \rightarrow \mathcal{C'}$ 
with $\c'$ a modular tensor category. A premodular
category is called \emph{modularizable}, if it admits a
modularization. 
\end{enumerate} 
\end{Definition} 

If a modularization exists, it is
unique up to equivalence of braided tensor categories.
It is known \cite[Corollary 3.5]{br} that a premodular
category over an algebraically closed field of characteristic
zero is modularizable, if and only if for all objects
$X \in M_{\c}$ one has $X \in T_{\c}$, $\theta_X = \id_X$ and 
$\dim X \in \mathbb{N}$. We thus obtain for the category $\mx$:

\begin{Proposition} \mbox{} \\
The premodular category $\mx$ is modularizable.
\end{Proposition}

\begin{Proof}
Corollary \ref{mt} implies for $p \in M_{\mx}$ 
that  $p \in T_{\mx}$ and $\theta_p = \id_p$. The assertion
follows by \cite[Corollary 3.5]{br}.
\end{Proof}

Proposition 2.3 of \cite{br} allows to detect
modularizations among dominant ribbon functors
$F: \c \rightarrow \c'$ between premodular categories:
it is sufficient to check that for any transparent object
$X \in M_{\c}$ the image $F(X)$ is trivial in the sense that
it is a finite direct sum of the tensor unit of $\c'$.

Let us investigate further the tensor subcategory of transparent objects:

\begin{Definition} \mbox{} \\
A premodular category $\c$ enriched over an algebraically closed
field $k$ is called \emph{Tannakian}, if there exists a modularization
of $\c$ that is equivalent to the category 
$\mathrm{vect}_f (k)$ of finite-dimensional $k$-vector spaces.

\end{Definition} 

We need the following facts proven in
\cite[Theorem 7.1]{De} and \cite[Theorem 2.11]{DM}:

\begin{Proposition}\label{tannakasch} \mbox{}\\
Let $\c$ be a premodular category over an algebraically closed field $k$ of characteristic
zero.
\begin{enumerate}
\item The category $\c$ is Tannakian, if and only if for all
simple objects  $X \in \Lambda_{\c}$ the twist equals the identity,
$\theta_X = \id_X$, and $\dim X \in \mathbb{N}$.

\item If $\c$ is Tannakian, it is equivalent as a tensor category
to the category of representations of a finite group $G$
on $k$-vector spaces.
\end{enumerate}
\end{Proposition}

\begin{Kor}\label{neutraltannakasch}\mbox{} \\
The full tensor subcategory $\mx^T$ of transparent objects of the premodular
tensor category $\mx$ is Tannakian.

\end{Kor}

\begin{Proof}
This follows immediately from Proposition  \ref{tannakasch} (i)
and Corollary \ref{mt} (ii).
\end{Proof}

We next determine explicitly the finite group $G=G(\x)$ describing the
Tannakian subcategory $\mx^T$. The action $\mu:\x_2 \times \x_1 \rightarrow \x_2$ that is part of the crossed module $\gm$
factorizes to an action of $\coker\,\partial$ on $\x_2$
which can be restricted to an action of $\coker\,\partial$
on $\ker\,\partial$:
\begin{equation}\label{hatmu}
        \begin{gathered}
                \ker\,\partial \times \coker\,\partial\\ 
                (k,Ig) 
        \end{gathered}
        \begin{gathered}
                \quad \rightarrow \quad \\
                \quad \mapsto \quad
        \end{gathered}
        \begin{aligned}
                &\ker\,\partial\\
                &k^{Ig} := k^g \,\,\, .
        \end{aligned}
\end{equation}  

Since the subgroup $\ker\,\partial$ is abelian, its irreducible 
characters form a group $(\ker\,\partial)^*$. 
We introduce the dual action

\begin{equation}
\begin{gathered}
                \hat \mu: (\ker\,\partial)^* \times \coker \,\partial\\ 
                (\chi,c) 
\end{gathered}
        \begin{gathered}
                \quad \rightarrow \quad \\
                \quad \mapsto \quad
        \end{gathered}
        \begin{aligned}
                &(\ker\partial)^*\\
                &\chi^c(k) := \chi(k^{c^{-1}})\,\, ,
        \end{aligned}
\end{equation} 
where we tacitly use the canonical identification
$(\ker\partial)^{**}\cong\ker\partial$.
We denote by $G(\x)$ the semi-direct product
\begin{equation}
G(\x):= (\ker\,\partial)^* \! \rtimes_{\hat \mu} \! 
(\coker\,\partial) \,\, .\label{G(X)}
\end{equation}

\begin{Proposition} \label{Gruppendarst} \mbox{} \\
The category $G(\x)\mbox{-}\mathrm{Rep}$ is equivalent, as a tensor category,
to the category $\mx^T$ of transparent $\x$-representations.
\end{Proposition}

\begin{Proof}
\begin{itemize}
\item
We construct the equivalence explicitly and define a functor
on objects as 
\begin{align}
F:G(\x)\mbox{-}\mathrm{Rep} \rightarrow \mathcal{M}{(\x)}\label{FunktorF}
\end{align}
which maps the $G(\x)$-representation $(V,\rho)$ 
to the triple $(V,P^{\rho},Q^{\rho})$ with
\begin{align*}
P^{\rho}(m) & := \left\{
\begin{array}{l  l}
                 \K \sum_{\chi \in K^*} \chi(m) \rho(\chi,I) & \textnormal{if } m \in K \equiv \ker \,\partial\\ 
                 0 &  \textnormal{else} 
        \end{array}
\right.\\[.2em]
Q^{\rho}(g) & := \rho(1,Ig).
\end{align*} 
Since a linear map commuting with the $G(\x)$-action commutes with
the action of $\x$ defined by $P^{\rho}$ and $Q^{\rho}$, we can
define $F$\ on morphisms as the identity so that the
functor $F$ is fully
faithful. To show that
the $\x$-representation $(V,P^{\rho},Q^{\rho})$ is transparent, consider any
$\x$-representation  $(W,P_W,Q_W)$ and compute the
braiding:
\begin{align*}
R_{V,W} & = \sum_{m \in \x_2} Q_W(\partial m) \otimes P^{\rho}(m) \circ \tau_{V,W} = \sum_{m \in K} Q_W(\partial m) \otimes \K \sum_{\chi \in K^*} \chi(m) \rho(\chi,I) \circ \tau_{V,W}\\[.2em]
& = Q_W(1) \otimes \sum_{\chi \in K^*} \delta(\chi,1) \rho(\chi,I) \circ \tau_{V,W} \qquad(\text{since} \; m \in \ker\,\partial)\\[.2em]
& = (\id_W \otimes \id_V) \circ \tau_{V,W} = \tau_{V,W}
\end{align*}
where $\tau_{V,W}:V\otimes W\to W\otimes V$ is the transposition map. 
Similarly, we find
\begin{align*}
R_{W,V} & = \sum_{m \in \x_2} Q^{\rho}(\partial m) \otimes P_W(m) \circ \tau_{W,V}
 = \rho(1,I) \otimes \sum_{m \inx} P_{W}(m) \circ \tau_{W,V}\\[.2em]
& = (\id_V \otimes \id_W) \circ \tau_{W,V} = \tau_{W,V}\,\, . 
\end{align*}

\item
We next show that $F$ is a strict tensor functor.
To check that the tensor unit of $G(\x)\mbox{-}\mathrm{Rep}$, the
trivial representation, is mapped to the tensor unit in
$\mx$, we remark that for $m \inx$ and $g \in \x_1$ we have
\begin{align*}
P^{\rho_{\mathbf{1}}}(m) & = \K \sum_{\chi \in K^*} \chi(m)\rho_{\mathbf{1}}(\chi,I) 
= \K\sum_{\chi \in K^*} \chi(m)\id_{\mathbb{C}} = \K\sum_{\chi \in K^*} m(\chi)\id_{\mathbb{C}} = \delta(m,1) \id_{\mathbb{C}},\\[.2em]
Q^{\rho_{\mathbf{1}}}(g) & = \rho_{\mathbf{1}}(1,Ig) = \id_{\mathbb{C}}.
\end{align*} 

Consider two objects $(V_1,\rho_1)$ and $(V_2, \rho_2)$ of
$G(\x)\mbox{-}\mathrm{Rep}$. To show that the tensor product of the
images under $F$ equals the image of the tensor product
$(V_1 \otimes V_2, P^{\rho_1 \otimes \rho_2}, Q^{\rho_1 \otimes \rho_2})$,
we remark
\begin{align*}
P_{V_{1}\otimes V_2}(m) & = \sum_{n \inx} P^{\rho_1}(n) \otimes P^{\rho_2}(n^{-1}m)\\[.2em]
& = \sum_{n \in K} \frac{1}{|K|^2}\sum_{\chi, \tilde \chi \in K^*} \chi(n) \tilde \chi(n^{-1}m) \rho_1(\chi,I) \otimes \rho_2(\tilde \chi,I)\\[.2em]
& = \K\sum_{\chi, \tilde \chi \in K^*} \delta(\chi,\tilde \chi) \chi(m)\rho_1(\chi,I) \otimes \rho_2(\tilde \chi, I)\\[.2em]
& = \K\sum_{\chi \in K^*} \chi(m)\rho_1(\chi,I) \otimes \rho_2(\chi,I) = P^{\rho_1 \otimes \rho_2}(m) \,\, , 
\end{align*}
where the third equality is the generalized orthogonality
relation for group characters \cite[Theorem 2.13]{is}.
The analogous identity for the action of $\x_1$ is straightforward.

\item
The functor $F$ being fully faithful, it suffices to show
that $F$ is essentially surjective to prove that it is an
equivalence of tensor categories.

Any transparent object is a direct sum of simple transparent
objects;\ hence we can restrict ourselves to simple
transparent objects.
They are all direct summands of the vacuum object
$\vk = (V_{\vk},P_{\vk},Q_{\vk})$ (Lemma \ref{vielfinvk}). From this, we conclude
that the linear map $P_{\vk}(m)$ is zero for
$m \notin K$ and that the automorphism $Q_{\vk}(g)$
is constant on the equivalence classes of the cokernel
$\coker\,\partial = \x_1 / I$.
Consider thus for a simple transparent object $(V,P,Q)$
\begin{equation}
\begin{array}{r c l}
        \rho:  G(\x) &  \rightarrow & \mathrm{Aut}(V)\\[.2em]
        (\chi,Ig) &     \mapsto & \rho(\chi,Ig) := \sum\limits_{k \in K} \chi^{-1}(k) Q(g)P(k).
\end{array}
\end{equation}
Direct computations show that this defines an action of the group $G(\x)$.

The image of the $G(\x)$-representation $(V,\rho)$ under $F$
is the $\x$-representation 
\begin{align*}
P^{\rho}(m) & =
\K \sum_{\chi \in K^*} \chi(m) \rho(\chi,I) = \K \sum_{k \in K, \chi \in K^*}  \chi(m) \chi^{-1}(k) Q(1)P(k)\\[.2em]
& = \K \sum_{k \in K^{**}, \chi \in K^*}  m(\chi)k(\chi^{-1}) P(k) = \sum_{k \in K^{**}} \delta(k,m) P(k) = P(m)& \text{if} \; m \in K\\
P^{\rho}(m) & = 0 = P(m) & \text{if} \; m \notin K\\
\intertext{and}
Q^{\rho}(g) & = \rho(1,Ig)
 = \sum_{k \in K} 1(k) Q(g)P(k) = Q(g) \,\, . 
\end{align*}
We conclude that $F$ is essentially surjective and thus an equivalence of tensor categories.
\end{itemize}
\end{Proof}

\section{The modularization  of $\mx$}

The vacuum object $\vk$ carries additional algebraic structure which
crucially enters in the modularization
of the premodular category $\mx$.

\begin{Definition}\label{Algebra} \mbox{} \\[-1.8em]
\begin{enumerate}
\item An \emph{algebra} in a (strict) tensor category $\c$ is
a triple consisting of an object $A\in\c$, a multiplication
morphism $m \in \hom(A \otimes A,A)$ and a unit
$\eta \in \hom(\teins, A)$ obeying the equations
$$ m \circ (m \otimes \id_A) = m \circ (\id_A \otimes m)
\qquad \text{ and }\qquad
m \circ (\eta \otimes \id_A) = \id_A = m \circ (\id_A \otimes \eta )\,\, .$$
A \emph{coalgebra} in $\c$ is defined analogously as a triple
consisting of an object $C$, a comultiplication morphism
$\Delta:C\to C\otimes C$ and a counit $\epsilon:C\to \teins$ obeying
coassociativity and counit equalities.
\item An algebra~$(A,m,\eta)$ in a braided tensor category $\c$ is
called (braided-)\emph{commutative}, if $m \circ R_{AA} = m$.

\item An algebra in a tensor category is called
\emph{haploid}, if it is simple as a left module over itself,
i.e.\ if $\dim_k \hom(\teins, A) = 1$.
\end{enumerate}
\end{Definition}

In the sequel we will see, that $\vk$ even carries the structure
of a special symmetric Frobenius algebra:

\begin{Definition}\label{Frobenius}\mbox{} \\
Let $\c$ be a (strict) tensor category. 
\begin{enumerate} 
\item A \emph{Frobenius algebra} in $\c$ is an object
with an algebra structure $(A,m,\eta)$ and a coalgebra
structure $(A, \Delta, \epsilon)$ such that 
$\Delta:A\to A\otimes A$  is a morphism
of $A$-bimodules:
\begin{align}
                        (\id_A \otimes m) \circ (\Delta \otimes \id_A) = \Delta \circ m = (m \otimes \id_A) \circ (\id_A \otimes \Delta).
                \end{align}
\item 
Suppose that the tensor category $\c$ is enriched over the category
of $k$-vector spaces where $k$ is a field.
A \emph{special algebra} in $\c$ is an object that is endowed
with an algebra and a coalgebra structure such that 
$$\epsilon \circ \eta  = \beta_1 \id_{\teins}
\quad\text{ and }\quad m \circ \Delta  = \beta_{A} \id_{A} $$
with invertible elements $\beta_1, \beta_A \in k^{\times}$.
    
\item 
Let $\c$ be a sovereign tensor category, i.e.\ a category
with left and right dualities that coincide as functors
from $\c$ to $\c^{\mathrm{opp}}$. A \emph{symmetric
algebra} in $\c$ is an algebra $(A,m,\eta)$ together with
a morphism $\epsilon \in \hom(A,\teins)$ such that the
two morphisms 

 $\Phi_1,\Phi_2: A \rightarrow A^{\vee}$
\begin{align}
\Phi_1 & := [(\epsilon \circ m) \otimes \id_{A^{\vee}}] \circ (\id_A \otimes b_A)\in\hom(A,A^{\vee})\\[.2em]
\Phi_2 & := [\id_{A^{\vee}} \otimes (\epsilon \circ m)] \circ (\tilde{b}_A \otimes \id_A)\in\hom(A,A^{\vee})
        \end{align}
are identical.

Here $b_A:\teins\to A\otimes A^\vee$ and
$\tilde b_A:\teins\to A^\vee\otimes A$ are the coevaluations
of the two dualities.
\end{enumerate}
\end{Definition}

Let $G$\ be a finite group and $k$ a field.
An important example of a symmetric special Frobenius 
algebra in the symmetric tensor category of $k[G]$-modules is
the algebra of functions $k(G)$ on $G$, the regular representation.

\begin{Lemma}\label{image}\mbox{} \\
The essential image of the regular representation of
$G(\x)$ under the functor $F$ is
the vacuum object $\vk$.
\end{Lemma}

\begin{Kor}\mbox{} \\
Since $k(G(\x))$ is a commutative symmetric Frobenius algebra,
the vacuum object $\vk$ carries a natural structure of a symmetric special Frobenius algebra
in $\mx^T$ and thus in $\mx$.
\end{Kor}
 
{\mbox{}\\\noindent \textbf{Proof of Lemma \ref{image}:{\hspace{0.7cm}
}\\}}
Consider the natural basis
$\{(\chi,c)\}_{(\chi,c) \in \kca}$ of $k(G(\x))$ 
of idempotents
\begin{align}
(\chi,c) \cdot (\tilde \chi, \tilde c) := \delta(\chi, \tilde \chi) \delta(c, \tilde c)(\tilde \chi, \tilde c)
\end{align}
in which the regular representation $\rho_R: G(\x)\to
\mathrm{Aut}(G(\x))$
is given by
\begin{align}
\rho_R(\chi,c) (\tilde \chi, \tilde c) := (\chi^{\tilde c} \tilde \chi, c \tilde c) \,\, .
\end{align}
It is convenient to perform a partial Fourier transform with respect to
$K$ to 
introduce also the basis
\begin{align}
(k,c) := \sum_{\chi \in K^*} \chi(k) (\chi,c)
\end{align}
of 

$k(G(\x))$
in which the multiplication is
$$ (k,c)\cdot (\tilde k,\tilde c)  
= \delta(c,\tilde c)(k \tilde k, \tilde c) \,\, . $$  
The regular algebra $k(G(\x))$ is mapped under the functor $F$ to the
triple\\ $(k(\kca),P^R,Q^R)$ with
\begin{align}
P^R(m) & = \left\{
\begin{array}{l  l}
                 \K \sum_{\chi \in K^*} \chi(m) \rho_R(\chi,I) 
                 & \textnormal{ if } m \in K\\ 
                 0 &  \textnormal{else} 
        \end{array}
\right.\\[.2em]
Q^R(g) & = \rho_R(1,Ig). 
\end{align} 
We compute the action of  $P^R$ and $Q^R$
on the basis $(k,c)_{(k,c) \in K \times C}$:

\begin{align*}
P^R(m) (\tk, I \tg) & = \K \sum_{\chi,\tilde \chi \in K^*} \chi(m) \tilde \chi(\tk)\rho_R(\chi,I)(\tilde\chi,I\tg)\\[.2em]
& = \K \sum_{\chi,\tilde \chi \in K^*} \chi(m^{\tg}) \tilde \chi(\tk)(\chi\tilde\chi,I\tg)
= \K \sum_{\chi, \chi' \in K^*} m^{\tg}(\chi) \tk(\chi^{-1}\chi')(\chi',I\tg)\\[.2em]
& = \sum_{\chi' \in K^*} \delta(m^{\tg},\tk) \chi'(\tk)(\chi',I\tg) = \delta(m^{\tg},\tk)(\tk,I\tg) & \text{if} \; m \in K\\
P^R(m) (\tk, I \tg) & = 0 = \delta(m^{\tg},\tk)(\tk,I\tg) & \text{if} \; m \notin K\\[.2em]
Q^R(g)(\tk, I \tg) &  = \sum_{\chi \in K^*} \chi(\tk)\rho_R(1,Ig)(\chi,I\tg) = \sum_{\chi \in K^*} \chi(\tk)(\chi,Ig\tg) = (\tk, I g\tg)
\,\, .
\end{align*}
Since this is precisely the action of $\x$ on $\vk$, we have 
proven the assertion.
\hfill$\Box$\\

Modules over the special symmetric Frobenius algebra $\vk$
crucially enter in the concrete construction 
\cite[Lemma 3.3]{br}
 of the
modularization:

\begin{Definition} \mbox{} \\
Let $A$\ be an algebra in a strict tensor category $\c$.
\begin{enumerate}
\item A (left) $A$-\emph{module} is a pair $(X,\rho_X)$ with $A\in\c$ and
$\rho_X \in \hom_{\c}(A \otimes X,X)$ such that 
$$\rho \circ (m \otimes \id_X) = \rho \circ ( \id_A \otimes \rho)\quad
\text{ and }\quad \rho \circ (\eta \otimes \id_X) = \id_X \,\, .
$$
\item A module $(X,\rho_X)$ over $A$ is called \emph{local} or
\emph{dyslectic} (\cite{Pa,KO,FFRS}),
if $\rho_X \circ R_{XA} \circ R_{AX} = \rho_X$.
\item A morphism of $A$-modules $(X,\rho_X)$ and $(Y,\rho_Y)$
is a morphism  $f \in \hom_{\c}(X,Y)$ such that
\begin{align}
                f \circ \rho_X = \rho_Y \circ (\id_A \otimes f)
                \,\, .
\end{align}
\item We denote by $A\mbox{-}\mod_{\c}$ the category of
$A$-modules in $\c$ and by $A\mbox{-}\mod_{\c}^{loc}$ 
the full subcategory of local $A$-modules.

\end{enumerate}
\end{Definition}

\begin{Remark} \mbox{} \\
Let $\c$\ be a braided tensor category and $A$\ be a commutative
algebra in $\c$.\ The following elementary facts from
commutative algebra are still valid in this setting:
\begin{enumerate}
\item Every left $A$-module $(M,\rho)$ has a 
structure of a right $A$-module with $(M,\rho\circ R_{M,A})$.

\item Let $M,N$ be two left $A$-modules. Then
$$\ M\otimes_A N:=\coker\,(\rho_M\circ R_{M,A}\otimes \id_N
- \id_M\otimes \rho_N)\ $$
endows the category $A\mbox{-}\mod_\c$ with the structure of a
tensor category.
\end{enumerate}
\end{Remark}

In fact, the modularization functor was constructed in 
\cite[Proposition 3.2]{br} as an induction functor for
a special commutative symmetric Frobenius algebra obtained
as the regular algebra in a Tannakian subcategory. We conclude

\begin{Proposition}\label{Induction}\mbox{} \\
The induction functor
\begin{align*}
\ind: \mx & \rightarrow \overline{\mx}:=\vk\mbox{-}\mod\\[.2em]
X & \mapsto (\vk \otimes X, \m \otimes \id_X) 
\end{align*}
is a modularization of  $\mx$.
\end{Proposition}

\begin{Proof}
By  \cite[Proposition 5.11]{FS},  the induction functor
is a tensor functor and by \cite[Proposition 5.17]{FS} it is
compatible with duality. Since $\vk$ is a special Frobenius
algebra, the induction functor is dominant by
\cite[Lemma 4.15]{FS}. From the explicit form of the braiding
given in \cite[Proposition 3.21]{FFRS} one deduces that
the induction functor  respects the braiding and is thus 
a ribbon functor. The category of modules over the regular algebra $k(G(\x))$ in $G(\x)\mbox{-}\mathrm{Rep}$
is equivalent to the category of $k$-vector spaces. Hence, for any transparent
object $X\in\c$, the induced module is isomorphic to a direct sum of $\vk$.
 Thus by 
\cite[Proposition 2.3]{br}, the induction functor is
a modularization.
\end{Proof}

\section{Explicit description of the modularization}

We now wish to describe the modularization $\overline{\mx}$
explicitly by showing that it is equivalent to the category
of representations of a crossed module $\bx$ with bijective boundary
map and thus to the representation category of an ordinary
Drinfeld double.
To this end, we consider 
\begin{align}
\bx := (I, \x_2/K, \overline{\mu}, \overline{\partial}) \label{xbar}
\end{align} 
with action
\begin{align*}
 \overline{\mu}: I \times \x_2/K & \rightarrow \x_2/K\\[.2em]
 (g,Km) & \mapsto K\mu(g,m) = Km^g
\end{align*}
and boundary map
\begin{align*}
\overline{\partial}: \x_2/K & \rightarrow I\\[.2em]
Km & \mapsto \partial(m)\,\, \, .
\end{align*} 
All maps are well-defined, since  $x \in K$ and $g \in I$ 
implies $x^g \in K$ and $\partial(x)= 1$. A direct computation
shows that this defines a crossed module;\ the bijectivity of
$\overline\partial$ is obvious.

\begin{Theorem}\label{Haupt} \mbox{} \\
The modularization of the representation category $\mx$ 
of a crossed module $\x$ is equivalent, as a ribbon category,
to the category of representations $\mbx$ of the crossed
module $\bx$.

\end{Theorem}

Our proof proceeds in two steps. We first introduce
the crossed module
\[\x' = (I, \x_2, \mu' = \mu|_{I \times \x_2}, \partial)\]
where we restrict to the image $I$ of $\partial$.
By abuse of notion, we denote
the boundary map of this crossed module again by
$\partial$;\ this map is surjective.
We denote by $\vk'$ the vacuum object of $\msx$ which is
again a commutative special symmetric Frobenius algebra.
We then construct a functor 
$$\vk \mbox{-}\mod_{\mx} \stackrel{F}{\longrightarrow} \vk' \mbox{-}\mod_{\msx}
\,\, . $$
Proposition \ref{Fverz} asserts that $F$\ is an equivalence
of abelian categories.

In a second step, we construct a functor
$$\vk' \mbox{-} \mod_{\msx} \stackrel{F'}{\longrightarrow} \mbx $$
and show in Proposition \ref{F'verz} that it provides
an equivalence of abelian categories as well. We finally
endow the two functors $F$ and $F'$ with the structure
of braided tensor functors and thus show that the categories
$\vk\mbox{-}\mod_{\mx}$ and $\mbx$ are equivalent as braided tensor categories.

This implies also that the categories are equivalent as
ribbon categories: any braided equivalence 
$G:\c \rightarrow \mathcal{D}$ of ribbon categories is
an equivalence of ribbon categories.
To see this, define on the image of $\c$ under $G$ a
new duality by $G(X)^* := G(X^{\vee})$. The new duality
is isomorphic to the duality in $\mathcal{D}$, thus
$G(X^{\vee}) \cong G(X)^{\vee}$. Since in any ribbon
category the twist can be expressed in terms of the dualities
and the braiding, the equivalence is also compatible with
the twist.

This concludes our argument that the categories
$\bmx$ and $\mbx$ are equivalent as ribbon categories.

We first construct a functor $F: \vk \mbox{-} \mod_{\mx}\rightarrow \vk' \mbox{-} \mod_{\msx}$ 
by restricting the group-action of $\x_1$ to the group-action of $I = \mathrm{Im} \partial$.\\

\noindent
\textbf{Construction of $F$}
\begin{itemize}
\item
To construct the functor $F$, we spell out the data
contained in an object of $\vk \mbox{-}\mod_{\mx}$.
Such an object consists of a $\x$-representation 
$(V,P_V,Q_V)$ and a $k$-linear map 
$\rho: V_{\vk} \otimes V \rightarrow V$ such that

\renewcommand{\labelenumii}{\roman{enumii})}
\begin{tabular}{l l }
\;\;\;\,(i)
$ \rho_V(x \otimes \delta_{Iy}, \rho_V(\tx \otimes \delta_{I\ty},v)) = \delta(Iy,I\ty)\rho_V(x\tx \otimes \delta_{I\ty},v)$&
($\vk$-action)\\[.3em]
\;\;\,(ii) 
$\rho_V(1 \otimes \sum_{Iy \in C} \delta_{Iy},v) = v$&(unitality of $\vk$-action)\\[.3em] 
\begin{tabular}{l}
(iii)
$\rho_V \circ P_{\vk V} = P_V \circ \rho_V$\\[.3em]
(iv)
$\rho_V \circ Q_{\vk V} = Q_V \circ \rho_V$
\end{tabular}&
($\rho_V$ is morphism in $\mx$)
\end{tabular}

\item
With the notation $V_{Iy} := 1 \otimes \delta_{Iy}.V$,
(i) and (ii) imply the decomposition of $V$ as
a direct sum of vector spaces
\[V = \bigoplus_{Iy \in C}V_{Iy} \,\, . \]
Similarly, we conclude that for every 
$x \in K$ the action $x \otimes \delta_{Iy}.\_$ is an
automorphism of vector spaces
\begin{align}
x \otimes \delta_{Iy}. V_{Iy} = V_{Iy} \,\, .\label{k-inv}
\end{align}
We next show that for all  $m \in \x_2$, $Iy \in C$, we have
\begin{align}
P_V(m)V_{Iy} \subset V_{Iy} \,\, . \label{x-inv}
\end{align}
Indeed,
\begin{align*}
P_{\vk V}(m)((x \otimes \delta_{Iy}) \otimes v) & = \sum_{n \in \mathcal{X}_2} P_{\vk}(n)(x \otimes \delta_{Iy}) \otimes P_V(n^{-1}m)v\\[.2em]
& = \sum_{n \in \mathcal{X}_2} \delta(n^y,x)(x \otimes \delta_{Iy}) \otimes P_V(n^{-1}m)v\\[.2em]
& = (x \otimes \delta_{Iy}) \otimes P_V((x^{y^{-1}})^{-1}m)v\\[.2em]
\end{align*}
and from (iii), we conclude
\[ P_V(m) (x \otimes \delta_{Iy}.v) = x \otimes \delta_{Iy}.(P_V((x^{y^{-1}})^{-1}m)v)
\,\,\, .\]

\item
From (iv) we conclude that for all  $Iy, I\ty \in C$,
we have vector space isomorphisms
$Q(\ty y^{-1}): V_{Iy} \rightarrow V_{I\ty}$
and that we have for all $h \in I$
\begin{align}
Q_V(h)V_{Iy} = V_{Iy} \,\, \, . \label{h-inv}
\end{align}
Indeed, we find with $g \in \x_2$
$$
Q_{\vk V}(g)((x \otimes \delta_{Iy}) \otimes v)  = Q_{\vk}(g)(x \otimes \delta_{Iy}) \otimes Q_V(g)v\\[.2em]
 = (x \otimes \delta_{Igy}) \otimes Q_V(g)v
$$
and thus by (iv)
\[ Q_V(g) (x \otimes \delta_{Iy}.v) = x \otimes \delta_{Igy}.(Q_V(g)v)
\,\, .\]

\item
From equations (\ref{k-inv}) - (\ref{h-inv}) we conclude that every
subvector space $V_{Iy}$ is invariant under the action
of  $x \in K$, $m \in \x_2$ and $h \in I$. In particular,
every vector space $V_{Iy}$ is a $\x'$-representation.
It becomes a $\vk'=\mathbb{C}[K]$-module by setting
$\rho'_V(x,v) := \rho_V(x \otimes \delta_{Iy}, v)$. 

All these $\vk'$-modules are isomorphic. We select the
$\vk'$-module $V_I$ as the image of the functor $F$:
$$
F(V,P_V,Q_V, \rho_V) := (V_{I},P_V,Q_V|_I, \rho_V(\_\otimes \delta_I,\_))
\,\, .
$$
On morphisms, we set
$$
F(\phi: V \rightarrow W) := \phi|_{V_I}: V_I \rightarrow W_I
\,\, \, . 
$$
Indeed, the image of the vector space $V_I$ under
$\phi:V \rightarrow W$ is contained in $W_I$, since 
$\phi$ commutes with the action, 
$\phi(V_I) = \phi(\rho_V(1 \otimes \delta_I,V)) = \rho_W(1 \otimes \delta_I,\phi(V)) \subset W_I$.
\end{itemize}

\begin{Proposition}\label{F"Aqu}\mbox{} \\
The functor $F$ presented in the above construction provides an equivalence of abelian
categories $\vk \mbox{-}\mod_{\mx}\simeq\vk' \mbox{-} \mod_{\msx}$.
\end{Proposition}

\begin{Proof}
We show that the functor is fully faithful and essentially surjective.
To show essential surjectivity, consider an object
$(W, P'_W,Q'_W,\rho'_W)$ in $\vk' \mbox{-}\mod_{\msx}$.

To find the preimage, we use induction from $I$ to $\x_1$: consider the object
$(V, P_V,Q_V,\rho_V)$ in $\vk' \mbox{-}\mod_{\msx}$
with 
\[V = \bigoplus_{Iy \in C}W_{Iy}\]
and action
\[ Q_V = \mathrm{Ind}_I^{\x_1}Q'_W: \x_1 \rightarrow 
\mathrm{End}(V) \,\, . \]
We introduce a $\x_2$-grading by 
\[P_V(m)w_{Iy} = (P_V(m)w)_{Iy}\] 
and the structure of a  $\vk$-module by
\[\rho_V(x \otimes \delta_{Iy},w_{I\ty}) := \delta(Iy,I\ty)(\rho'_W(x,w))_{Iy}
\,\, \, . \]
A straightforward calculation shows that the image
of this object under $F$\ is $(W, P'_W,Q'_W,\rho'_W)$.

To show that $F$ is fully faithful, we note that a morphism
$\phi: V \rightarrow W$ from $(V,P_V,Q_V,\rho_V)$ to 
$(W,P_W,Q_W,\rho_W)$ is determined by its restriction to
$V_I$, since for any $v\in V$, we have
\begin{align*}
\phi(v) & = \sum_{Iy \in C} \phi(1 \otimes \delta_{Iy}.v)
 = \sum_{Iy \in C} \phi( Q(y) Q(y^{-1})1 \otimes \delta_{Iy}.v)\\[.2em]
& = \sum_{Iy \in C} Q(y) \phi\underbrace{(1 \otimes \delta_{I}. 
Q(y^{-1})v)}_{\in V_I} \,\,\, .
\end{align*} 
\end{Proof} 

We next construct an equivalence  $F':\vk' \mbox{-}\mod_{\msx} \rightarrow \mbx$.
The idea is to take coinvariants with respect to the action of the
kernel $K:=\ker\,\partial$.\\

\noindent
\textbf{Construction of $F'$}
\begin{itemize}
\item
To construct an equivalence
 $F':\vk' \mbox{-} \mod_{\msx} \rightarrow \mbx$ we spell out
 explicitly the data contained in an object
$(W, P_W, Q_W, \rho_W)$ of $\vk'\mbox{-}\mod_{\msx}$:
here $(W, P_W, Q_W)$ is an object of $\msx$ and 
$\rho_W: \mathbb{C}[K] \otimes W \rightarrow W$ 
is a $k$-linear map such that
\begin{enumerate}
\item $\rho_W(x,\rho_W(\tx,w)) = \rho(x \tx,w)$
\item $\rho_W(1,w) = w$
\item $\rho_W \circ P_{\vk' W} = P_W \circ \rho_W$
\item $\rho_W \circ Q_{\vk' W} = Q_W \circ \rho_W$
\end{enumerate}
We introduce the shorthand notation  $x.w$ for $\rho_W(x,w)$.

\item
From the first two axioms we conclude that for all
$x\in K$, the action is an isomorphism. As a
$\x'$-module, we decompose
\[ W = \bigoplus_{m \in \x_2} W_m
\quad\text{ with }\quad W_m := P(m)W \,\, .
\]
From (iii) we conclude as in the construction of $F$
\begin{align}
x.P(m)w = P(xn)x.w \,\,\, . \label{xp}
\end{align}
Thus the action of $x$ implies for
$Km = Kn$ in $\x_2 /K$ the isomorphy of vector spaces
$W_m\cong W_n$. Again as in the construction of $F$,
we conclude
\begin{align}
x.(Q(g)w) = Q(g)(x.w) \quad\text{for all} \quad g\in \x_1
\,\, . \label{xq}
\end{align}
\item
We now define $F'$ by taking coinvariants with respect
to the action of $K=\ker\,\partial$. On objects, we have
\[F'(W,P_W,Q_W,\rho_W) := (W_K,P_{W_K},Q_{W_K})\]
with
\[ W_K := W/(x.w - w| x \in K, w \in W)\]
\[P_{W_K}(Km)\overline{w} := \overline{P_W(m)w}\]
\[Q_{W_K}(h)\overline{w} := \overline{Q_W(h)w} \,\, . \]
From (\ref{xp}) and (\ref{xq}) we deduce that the maps
$P_{W_K}$ and $Q_{W_K}$ are well-defined.
On morphisms, we consider the restriction
\[F(f: W \rightarrow V) := (f_K: W_K \rightarrow V_K) \]
and obtain a $k$-linear functor $F'$.
\end{itemize}

\begin{Proposition}\label{F'"Aqu} \mbox{} \\
The functor $F'$ presented in the above construction provides an equivalence of abelian
categories  $\vk' \mbox{-} \mod_{\msx}\simeq\mbx$. 
\end{Proposition}

\begin{Proof}
To show that $F'$ is essentially surjective,
we construct for $(V, P_V, Q_V)\in\mbx$ an object 
$(W,P_W, Q_W, \rho_W)\in\vk'\mbox{-}\mod_{\mx}$ with
\[W := \bigoplus_{x \in K} W_x \,\, , \]
where $W_x \cong V$ as a $\x_1$-representation for all 
$x \in K$. On $W$ we define an action of $K$ by 
\begin{align*}
x.w_{\tx} := w_{x \tx} \,\, \, . 
\end{align*} 
To define an action of $\x_2$, choose representatives
$(Km)$ and set
\begin{align*}
P_W(xm)w_{\tx} := \delta(x,\tx)P_{W_x}(Km)w_x \,\, . 
\end{align*}
The image under $F$\ of the object constructed is isomorphic
to $(V, P_V, Q_V)$, showing essential surjectivity.

As in the proof of proposition \ref{F"Aqu}, we conclude that
a morphism $\phi: V \rightarrow W$ is uniquely determined 
by the map $\phi_K$ induced on coinvariants.
\end{Proof}

It remains to endow the functors with more structure.

\begin{Proposition}\label{Fverz} \mbox{} \\
Consider the morphism
\begin{equation*}
        \begin{gathered}
                \varphi_0:\\
                {}
        \end{gathered}
        \begin{gathered}
                \mathbb{C}[K]\\
                x
        \end{gathered}
        \begin{gathered}
                \rightarrow\\
                \mapsto
        \end{gathered}
        \begin{aligned}
                &\mathbb{C}[K] \otimes \mathbb{C}(\delta_I)\\
                &x \otimes \delta_I 
        \end{aligned}
\end{equation*}
and for all objects $V,W$ of $\vk'\mbox{-}\mod_{\msx}$ the morphisms: 
\begin{equation*}
        \begin{gathered}
                \varphi_2(V,W):\\
                {}
        \end{gathered}
        \begin{gathered}
                F(V) \otimes_{\vk'} F(W)\\
                v_I \otimes_{\vk'} w_I 
        \end{gathered}
        \begin{gathered}
                \rightarrow\\
                \mapsto
        \end{gathered}
        \begin{aligned}
                &F(V \otimes_{\vk} W)\\
                &v_I \otimes_{\vk} w_I \,\,\, .
        \end{aligned}
\end{equation*}
These morphisms endow the functor
$F: \vk \mbox{}-\mathrm{Mod}_{\mx} \rightarrow \vk' \mbox{-}\mathrm{Mod}_{\msx}$ 
with the structure of a braided tensor functor.
\end{Proposition}

\begin{Proof}
Bijectivity of $\varphi_0$ is obvious. To check
bijectivity of $\varphi_2(V,W)$, we note that the equation\\
$v \otimes_{\vk} w = \delta_I.(v \otimes_{\vk} w) = 
\delta_I.v \otimes_{\vk} \delta_I.w$
for $v \otimes_{\vk} w \in F(V \otimes W)$ 
implies
$v \otimes_{\vk} w = \varphi_2(\delta_I.v \otimes_{\vk'} \delta_I.w)$
and hence surjectivity. On the other hand, 
$\varphi_2(V,W)(v_I \otimes_{\vk'} w_I) = 0$
implies $v_I = 0$ and $w_I = 0$ and thus injectivity
of $\varphi_2$. The verification that 
$(F,\varphi_0, \varphi_2)$ is a tensor functor is routine.

To show that the tensor functor $(F,\varphi_0, \varphi_2)$
is braided, we have to check that for any pair of objects
$(V,P_V,Q_V), (W,P_W,Q_W)$ the diagram
\[\begin{CD}
F(V) \otimes_{\vk'} F(W) @>\varphi_2 >> F(V\otimes_{\vk} W)\\
@VV R_{F(V),F(W)} V @VV F(R_{V,W}) V\\
F(W) \otimes_{\vk'} F(V) @>\varphi_2>> F(W \otimes_{\vk} V)
\end{CD}\]
commutes; indeed, 
\begin{align*}
F(R_{V,W}) \circ \varphi_2(V,W)(v \otimes_{\vk'} w) & = F(R_{V,W}) (v \otimes_{\vk} w)\\[.2em]
& = \sum_{n \inx} Q_W(\partial' m) w \otimes_{\vk} P_V(m)v \\[.2em]
& = \varphi_2 \left( \sum_{n \inx} Q_W(\partial' m) w \otimes_{\vk'} P_V(m)v \right) \\[.2em]
& = \varphi_2 \circ R_{F(V),F(W)}(v \otimes_{\vk'} w)\,\, .
\end{align*}
Hence $(F, \varphi_0, \varphi_2)$ is a braided tensor functor.

\end{Proof}

\begin{Proposition}\label{F'verz}
Consider the morphisms
\begin{equation*}
        \begin{gathered}
                \varphi'_0:\\
                {}
        \end{gathered}
        \begin{gathered}
                \mathbb{C}\\
                \lambda
        \end{gathered}
        \begin{gathered}
                \rightarrow\\
                \mapsto
        \end{gathered}
        \begin{aligned}
                &(\mathbb{C}[K])_K\\
                &\lambda \overline{x} & x \in K
        \end{aligned}
\end{equation*}
and for all objects $V,W$ of $\mbx$
\begin{equation*}
        \begin{gathered}
                \varphi'_2(V,W):\\
                {}
        \end{gathered}
        \begin{gathered}
                V_K \otimes W_K \\
                \overline{v} \otimes \overline{w} 
        \end{gathered}
        \begin{gathered}
                \rightarrow\\
                \mapsto
        \end{gathered}
        \begin{aligned}
                &(V \otimes_{\vk'} W)_K\\
                &\overline{v \otimes_{\vk'} w}\,\, .
        \end{aligned}
\end{equation*}
These morphisms endow the functor 
$F': \vk'\mbox{-}\mod_{\msx} \rightarrow \mbx$ with the structure
of a braided tensor functor.
\end{Proposition}

\begin{Proof}
We first remark that
$\varphi'_0$ is well-defined, since for $x, x' \in K$ 
we have $\overline{x} = \overline{x'}$
in $(\mathbb{C}[K])_K$. The bijectivity of $\varphi_0$
is immediate from $\dim_k(\mathbb{C}[K])_K = 1$ and 
$\ker \,\varphi'_0 = {0}$.

To check that also $\varphi_2'(V,W)$ is well-defined, we first
remark that the action of $x \in K$ on
$v \otimes_{\vk'} w\in V \otimes_{\vk'} W$ reads
\[x.(v \otimes_{\vk'} w) = x.v \otimes_{\vk'} w 
= v \otimes_{\vk'} x.w\,\, . \]
Now take $v,v'\in V$ and $w,w'\in W$ such that 
\[\overline{v} \otimes \overline{w} = \overline{v'} \otimes \overline{w'}
\,\, .\]
Then we can find  $x, \tx \in K$ such that  $v' = x.v$ and 
$w' =  \tx.w$ and we have by the proceeding remark
\[\overline{v \otimes_{\vk'} w} = \overline{x \tx.(v \otimes_{\vk'} w)} = \overline{x.v \otimes_{\vk'} \tx.w}\]
and thus
\[\varphi_2'(\overline{v}\otimes \overline{w}) = \varphi_2'(\overline{v'}\otimes \overline{w'})\,\,.\]
An inverse of $\varphi_2'$ can be given directly by
\[\varphi_2^{'-1}(\overline{v \otimes_{\vk'} w}) = \overline{v} \otimes \overline{w}
\,\, .\]

One checks by direct computations that 
$(F', \varphi_0', \varphi_2')$ is a tensor functor. Finally,
$(F', \varphi_0', \varphi_2')$ is braided, since we have
\begin{align*}
F'(R_{VW}) \circ \varphi_2'(\overline{v} \otimes \overline{w}) & = F'(R_{VW})(\overline{v \otimes_{\vk'} w})\\[.2em]
& = \sum_{Km \in \x_2/K} \overline{Q(\bar\partial Km)w \otimes_{\vk'} P(Km)v}\\[.2em]
& = \varphi_2' \left( \sum_{Km \in \x_2/K} \overline{Q(K \bar\partial m)w} \otimes \overline{P(Km)v} \right) \\[.2em]
& = \varphi_2' \circ R_{V_K W_K}(\overline{v} \otimes \overline{w})
\,\,.
\end{align*}
\end{Proof}
 
\subsubsection*{Acknowledgements}
We thank J\"urgen Fuchs, Thomas Nikolaus and Ingo Runkel for helpful discussions.
The authors are partially supported by the DFG Priority Program 1388
``Representation Theory''.


\begin{thebibliography}{---------}
\bibitem[BK]{BK}
B. Bakalov, A. Kirilov,
\emph{Lectures on tensor categories and modular functors}, 
American Mathematical Society, Providence, 2000

\bibitem[Ba]{bantay}
P. Bantay, 
\emph{Characters of crossed modules and premodular categories},\\
London Mathematical Society Lecture Note Series 372 (2010)
{\tt math.QA/0512542}

\bibitem[Br]{br}
A. Brugui\`eres, \emph{Cat\'egories pr\'emodulaires, modularisations et invariants des vari\'et\'es de dimension 3}, Math. Ann. 316 (2000) 215-236

\bibitem[De]{De}
P. Deligne, \emph{Cat\'egories tannakiennes} in \emph{The Grothendieck Festschrift Volume II}, 
Progr. Math., 87, Birkh\"auser Boston, Boston, MA, 1990, 111-195

\bibitem[DM]{DM}
P. Deligne, J.S. Milne, \emph{Tannakian categories} in \emph{Hodge cycles, motives, and Shimura varieties},
Lecture Notes in Mathematics 900, Springer, Berlin/Heidelberg, 1982,
101-228

\bibitem[DVVV]{DVVV}
R.\ Dijkgraaf, C.\ Vafa, E.\ Verlinde, and H.\ Verlinde,
\emph{The operator algebra of orbifold models},
Commun.\ Math.\ Phys.\ 123 (1989) 485-526

\bibitem[FFRS]{FFRS}
J. Fr\"ohlich, J. Fuchs, I. Runkel, C. Schweigert,
\emph{Correspondences of ribbon categories},
Adv. Math.  199  (2006) 192-329,
{\tt math.CT/0309465}

\bibitem[FS]{FS}
J. Fuchs, C. Schweigert,
\emph{Category theory for conformal boundary conditions},
Fields Institute Communications 39 (2003) 25-71,
{\tt math.CT/0106050}

\bibitem[Is]{is}
I.M. Isaacs,
\emph{Character theory of finite groups},
Academic Press, New York, 1977

\bibitem[Ka]{Ka}
C. Kassel, 
\emph{Quantum groups}, 
Springer, Graduate Texts in Mathematics 155, New York, 1995

\bibitem[KO]{KO}
A. Kirillov, V. Ostrik, 
\emph{On $q$-analog of McKay correspondence and ADE classification of $\widehat{\frak{sl}}_2$ conformal field theories},
Adv. Math. 171 (2002) 183-227, {\tt math.QA/0101219}

\bibitem[M1]{M1}
M.\ M\"uger,
\emph{Galois theory for braided tensor categories and the modular closure},
Adv.\ Math.\ 150 (2000) 151-201, {\tt math.CT/9812040}


\bibitem[N]{N}
D.\;Naidu, 
\emph{Crossed pointed categories and their equivariantizations}, 
Pacific J. Math. 247 (2010), no. 2, 477–-496. 


\bibitem[Pa]{Pa}
B. Pareigis, 
\emph{On braiding and dyslexia}, 
J.\ Algebra 171 (1995) 413�425

\bibitem[RT]{RT}
N.Yu.\ Reshetikhin and V.G.\ Turaev,
\emph{Invariants of $3$-manifolds via link polynomials and quantum groups}
Inv.\ Math.\ 103 (1991) 547-597

\bibitem[Tu]{Tu}
V.G. Turaev,
\emph{Quantum invariants of knots and 3-manifolds}, 
W. de Gruyter,
Berlin, 1994

\end{thebibliography}
\end{document}